\documentclass[a4paper,11pt]{article} 

\usepackage[utf8]{inputenc} % Use UTF-8 encoding
\usepackage{microtype} % Slightly tweak font spacing for aesthetics
\usepackage[english]{babel} % Language hyphenation and typographical rules
\usepackage{amsthm, amsmath, amssymb} % Mathematical typesetting
\usepackage{float} % Improved interface for floating objects
\usepackage[final, colorlinks = true, 
linkcolor = black, 
citecolor = black]{hyperref} % For hyperlinks in the PDF
\usepackage{caption, subcaption}
\usepackage{graphicx, multicol} % Enhanced support for graphics
\usepackage{xcolor} % Driver-independent color extensions
\usepackage{rotating} % Rotation tools
\usepackage{censor} % Facilities for controlling restricted text
\usepackage{matlab-prettifier} % Environment for specifying algorithms in a natural way
\usepackage{booktabs} % Enhances quality of tables
\usepackage[yyyymmdd]{datetime} % Uses YEAR-MONTH-DAY format for dates
\usepackage{fancyhdr} % Headers and footers
\usepackage{pgfplots}
\pgfplotsset{compat=1.15}
\usepackage{hyperref}
\usepackage{tikz-cd,mathtools}

\usetikzlibrary{arrows}
\usetikzlibrary{braids}

\usepackage{wrapfig}
\usepackage{pdfpages}
\usepackage{mwe}

\usepackage{caption}

\DeclareMathOperator{\R}{\mathbb R}

\DeclareMathOperator{\N}{\mathbb N}
\DeclareMathOperator{\Q}{\mathbb Q}
\DeclareMathOperator{\Z}{\mathbb Z}

\DeclareMathOperator{\rank}{rank}
\DeclareMathOperator{\im}{Im}

\DeclareMathOperator{\Ker}{Ker}
\DeclareMathOperator{\Span}{Span}

\theoremstyle{plain}
\newtheorem{theorem}{Theorem}
\newtheorem{lemma}[theorem]{Lemma}
\newtheorem{proposition}[theorem]{Proposition}
\newtheorem{corollary}[theorem]{Corollary}
\newtheorem{conjecture}[theorem]{Conjecture}

\newtheorem{question}[theorem]{Question}
\newtheorem{fact}[theorem]{Fact}
\theoremstyle{definition}
\newtheorem{remark}[theorem]{Remark}

\newtheorem{construction}[theorem]{Construction}
\newtheorem{definition}[theorem]{Definition}
\newtheorem{notazione}[theorem]{Notation}
\newtheorem{example}[theorem]{Example}
\newtheorem{exercise}[]{Exercise}
\newtheorem{problem}[theorem]{Problem}

\newtheorem{vuoto}[theorem]{}

\numberwithin{theorem}{section}

\newcommand{\bt}{\begin{theorem}}
\newcommand{\et}{\end{theorem}}

\newcommand{\bv}{\begin{vuoto}}
\newcommand{\ev}{\end{vuoto}}

\newcommand{\bl}{\begin{lemma}}
\newcommand{\el}{\end{lemma}}

\newcommand{\bd}{\begin{definition}}
\newcommand{\ed}{\end{definition}}

\newcommand{\beq}{\begin{equation}}
\newcommand{\eeq}{\end{equation}}

\newcommand{\bexa}{\begin{example}}
\newcommand{\eexa}{\end{example}}

\newcommand{\bexe}{\begin{exercise}}
\newcommand{\eexe}{\end{exercise}}

\newcommand{\bfact}{\begin{fact}}
\newcommand{\efact}{\end{fact}}

\newcommand{\bprop}{\begin{proposition}}
\newcommand{\eprop}{\end{proposition}}

\newcommand{\bp}{\begin{proof}}
\newcommand{\ep}{\end{proof}}

\newcommand{\bc}{\begin{corollary}}
\newcommand{\ec}{\end{corollary}}

\newcommand{\bq}{\begin{question}}
\newcommand{\eq}{\end{question}}

\newcommand{\bcong}{\begin{conjecture}}
\newcommand{\econg}{\end{conjecture}}

\newcommand{\bproblem}{\begin{problem}}
\newcommand{\eproblem}{\end{problem}}

\newcommand{\bs}{\begin{proof}[Proof.]}
\newcommand{\es}{\end{proof}}

\newcommand{\br}{\begin{remark}}
\newcommand{\er}{\end{remark}}

\newcommand{\bn}{\begin{notazione}}
\newcommand{\en}{\end{notazione}}

\begin{document}

\title{The Thurston norm of graph manifolds}
\author{Alessandro V. Cigna}
\date{\today}
\maketitle

\begin{abstract}
    \noindent The Thurston norm of a closed oriented graph manifold is a sum of absolute values of linear functionals, and either each or none of the top-dimensional faces of its unit ball are fibered. We show that, conversely, every norm that can be written as a sum of absolute values of linear functionals with rational coefficients is the nonvanishing Thurston norm of some graph manifold, with respect to a rational basis on its second real homology. Moreover, we can choose such graph manifold either to fiber over the circle or not. In particular, every symmetric polygon with rational vertices is the unit polygon of the nonvanishing Thurston norm of a graph manifold fibering over the circle. In dimension $\ge 3$ many symmetric polyhedra with rational vertices are not realizable as nonvanishing Thurston norm ball of any graph manifold. However, given such a polyhedron, we show that there is always a graph manifold whose nonvanishing Thurston norm ball induces a finer partition into cones over the faces.
\end{abstract}

\section{Introduction}
The Thurston norm \cite{norm} is a powerful tool for analyzing the topology of a 3-manifold by studying embedded surfaces sitting inside the $3$-manifold.
Given a compact orientable 3-manifold $M$, the Thurston norm is a seminorm $x$ defined on the vector space $H_2(M,\partial M;\R)$. If $S$ is some properly embedded oriented surface in $M$, then $-x([S])$ is the maximal Euler characteristic among the properly embedded oriented surfaces homologous to $S$, computed after discarding any component being a sphere, a disc, a torus or an annulus.  For instance, in a knot complement, the Thurston norm detects the minimal genus of the knot. 

The unit ball $B_x$ of $x$ is a finite (possibly unbounded) convex polyhedron (cf. \cite{norm}), thus the whole information carried by $x$ is determined by a finite amount of data, namely the combinatorics of the faces of $B_x$. 

\bigskip

Unfortunately, the Thurston norm has been described just in some special families of manifolds (see for example \cite{norm, pretzels, os, tunnel, chen, cigna}). In \cite{cigna}, we focused on the manifolds obtained as $2$-bridge link complements. There, the Thurston norm is a seminorm on the plane and happens to be always quite simple: the unit polygon has at most $8$ vertices, lying along the axes and the principal bisectors. However, by means of this result, we showed that the Thurston norms obtainable on the plane can have unit balls with an arbitrary number of sides, even if we restrict to complements in $S^3$ of links with two unknotted components. 

In the present paper, we move our attention to a different family of $3$-manifolds, namely graph manifolds. In many cases, the second homology of a graph manifold has dimension $\ge 3$, a setting where even less is known about the possible shapes of Thurston unit polyhedra. The advantage of graph manifolds is that their surfaces can be easily understood. Indeed, a graph manifold decomposes into Seifert fibered pieces, and an incompressible and boundary-incompressible surface in a Seifert fibered manifold can be isotoped to become either \emph{horizontal} or \emph{vertical}. Moreover, a useful result by Neumann \cite{neumann} helps us to understand when surfaces inside the different pieces can patch together to give a properly embedded surface in the ambient graph manifold.

Let $M$ be a closed oriented graph manifold. In Section \ref{sec: norm of a graph manifold} we show that the Thurston norm on $H_2(M;\R)$ can be written as a sum of absolute values of linear functionals, and either no or all top-dimensional faces of the Thurston ball are fibered (see Lemma \ref{lemma: norm} and Proposition \ref{prop: fibering}). The main result of this article is the following.

\bigskip

\noindent\textbf{Theorem \ref{thm: every sum is graph}.} \textit{Let $||\cdot||$ be a norm on $\R^d$ that can be written as $$||\cdot ||:=\sum_{i=1}^n |\langle \beta_i,\cdot\rangle|$$ for some $\beta_1,...,\beta_n\in \Q^d$, and let $g_1,...,g_n$ be nonnegative integers. There exist a closed oriented good graph manifold $M$ with base surfaces of genera $g_1,...,g_n$ and a rational basis on $H_2(M;\R)$, such that the nonvanishing Thurston norm on $H_2^{nv}(M;\R)$ in the induced basis coincides with $||\cdot||$. Moreover, we can ask for either all or none of the top-dimensional Thurston cones of $H_2(M;\R)$ to be fibered.}

\bigskip

In Theorem \ref{thm: every sum is graph}, \emph{good} graph manifold means that its Seifert fibered JSJ pieces all have orientable base orbifold, whereas the \emph{nonvanishing homology $H_2^{nv}(M;\R)$} and the \emph{nonvanishing Thurston norm} refer to the quotient of $H_2(M;\R)$ by the null-space of the Thurston norm and the induced norm respectively. Notice that, in the case of good graph manifolds, we can also compute the dimension of the null-space of the Thurston norm.

\bprop[See Proposition \ref{prop: dimension}] Suppose $M$ is a closed oriented good graph manifold, and let $\Gamma$ be the decomposition graph of $M$. For $i=1,...,n$, let $g_i$ be the genus of the base orbifold of the piece $M_i$. The null-space of the Thurston norm on $H_2(M;\R)$ has dimension $|\Gamma|-\chi(\Gamma)+2\sum_{i=1}^n g_i$, where $|\Gamma|$ is the number of connected components of $\Gamma$ and $\chi(\Gamma)$ is its Euler characteristic. 
\eprop

Then, the natural question is how ``generic" the norms of Theorem \ref{thm: every sum is graph} are among the norms with polyhedral unit ball. In two dimensions the answer is optimal, i.e. every norm with polygonal unit ball is a sum of absolute values of linear functionals. Thus, we recover the following result.

\bigskip

\noindent\textbf{Corollary \ref{cor: planar norms}.} \textit{If $||\cdot ||$ is a norm on $\R^2$ whose unit ball is a polygon with rational vertices, there exist a closed oriented good graph manifold $M$ and a rational basis on $H_2(M;\R)$, such that the nonvanishing Thurston norm on $H_2^{nv}(M;\R)$ coincides with $||\cdot ||$ in the induced basis. Moreover, we can ask for either all or none of the top-dimensional Thurston cones of $H_2(M;\R)$ to be fibered.}

\bigskip

Observe that, in \cite{norm}, Thurston has shown a result that is stronger than Corollary \ref{cor: planar norms} as it just uses integral bases and gives ways to realise the polygons directly on $H_2(M;\R)$, in such a way that any prescribed collection of edges corresponds exactly with fibered cones. At the same time, the nice feature of Corollary \ref{cor: planar norms} is that we can keep track of more structure on the manifolds used to realise the given planar norms, namely their being good graph manifolds together with whether they fiber over the circle or not.

\bigskip

In higher dimension, there are many norms whose unit ball is a polyhedron that cannot be written in the form of Theorem \ref{thm: every sum is graph}. If a norm on $\R^n$ has polyhedral unit ball $P$, then we can study the cellular decomposition of the sphere $S^{n-1}$ induced by projecting radially the cellular decomposition of $\partial P$ onto it. In Section \ref{sec: cellular decompositions}, we define a notion of \emph{completeness} for a cellular decomposition of $S^{n-1}$ (see Definition \ref{def: completeness}). Then a \emph{complete polyhedron} is a polyhedron inducing a complete cellular decomposition of the sphere. Many polyhedra are not complete, and the cube $[-1,1]^n$ is one such example for $n\ge 3$. By the following proposition, such polyhedra are not realizable as (nonvanishing) Thurston balls of any closed oriented graph manifold.

\bigskip

\noindent\textbf{Proposition \ref{prop: completeness}} \textit{Let $\beta_1,...,\beta_k\in \R^n-\{0\}$ be a set of generators of $\R^n$. Let $||\cdot||$ be the norm $$||v||:=\sum_{i=1}^k|\langle \beta_i,v\rangle|.$$ Then, $||\cdot||$ has polyhedral unit ball, the linear cellular decomposition of $S^{n-1}$ induced by $||\cdot||$ is complete and its $1$-codimensional skeleton is} $$S^{n-1}\cap\left(\bigcup_{i=1}^k \beta_i^{\bot}\right).$$

\bigskip

Also, there are examples of $3$-manifolds whose Thurston ball is not complete, so the (nonvanishing) norms of graph manifolds do not exhaust all the possible norms obtainable as (nonvanishing) Thurston norm of some $3$-manifold (see Example \ref{exa: 3-chain}).

Finally, not every complete polyhedron is the nonvanishing Thurston ball of a graph manifold (see Remark \ref{rem: icosaedro}). However, there is a \emph{completion} construction which makes it possible to state that, given a norm with polyhedral unit ball, there is always a ``more complicated" norm which arises as Thurston norm of a graph manifold (see the discussion following Corollary \ref{cor: completion}).

\bigskip

\noindent\textbf{Corollary \ref{cor: completion2}} \textit{ Let $P\subset \R^n$ be a polyhedron with rational vertices which is symmetric through the origin. There is a closed oriented good graph manifold $M$, such that there is a rational basis in which $H_2^{nv}(M;\R)\cong \R^n$ and the partition of $H_2^{nv}(M;\R)$ into cones over the faces of the nonvanishing Thurston unit ball is a refinement of the partition of $\R^n$ into cones over the faces of $\partial P$. Moreover, we can ask for either all or none of the top-dimensional Thurston cones of $H_2(M;\R)$ to be fibered.}

\bigskip

\noindent \textbf{Structure of the paper.}  In Section \ref{sec: preliminaries} we recall the definitions and basic properties of the main objects on focus: the Thurston norm, Seifert fibered manifolds, and graph manifolds. Along the way, we introduce a result by Neumann \cite{neumann} that facilitates our understanding of the surfaces sitting inside graph manifolds. In Section \ref{sec: norm of a graph manifold} we compute the norm of a graph manifold. In Section \ref{sec: characterization}, a characterization of all norms that can be obtained as nonvanishing Thurston norm of graph manifolds is achieved (see Theorem \ref{thm: every sum is graph}). Finally, in Section \ref{sec: cellular decompositions} we focus on the polyhedra arising as unit balls of Thurston norms on graph manifolds. This last section aims to understand the geometry of such polyhedra and their properties better. In this way, we will have an idea of how ``generic" such polyhedra are, inside the class of all polyhedra or the class of those polyhedra arising as Thurston unit balls of some $3$-manifold. 

\bigskip

\noindent\textbf{Aknowledgements.} I would like to thank my PhD advisor Mehdi Yazdi for useful discussions on the
topic and for his supervision. Many thanks also to the anonymous referee for their thoughtful comments on the initial version of this paper.

\section{Preliminaries}\label{sec: preliminaries}
\subsection{(Nonvanishing) Thurston norm}
We now recall the basic definitions and properties of the Thurston norm. For further details, we refer to \cite{norm}. From now on, $M$ will always be a $3$-dimensional compact orientable smooth manifold.

\bigskip

\noindent Every properly embedded oriented surface $S\subset M$ represents a class $[S]\in H_2(M,\partial M; \Z)$. Conversely, every class in $H_2(M,\partial M;\Z)$ is represented by a properly embedded oriented surface. 

A quite natural question is to ask whether the sum in $H_2(M,\partial M;\Z)$ corresponds with some operation between surfaces in $M$. The answer is affirmative: if $[S_1]=\alpha_1$ and $[S_2]=\alpha_2$, then $\alpha_1+\alpha_2$ is represented by the \emph{oriented cut-and-paste} of $S_1$ and $S_2$. Here is the definition.

\begin{figure}[ht]
    \centering
    \includegraphics[width=0.85\textwidth]{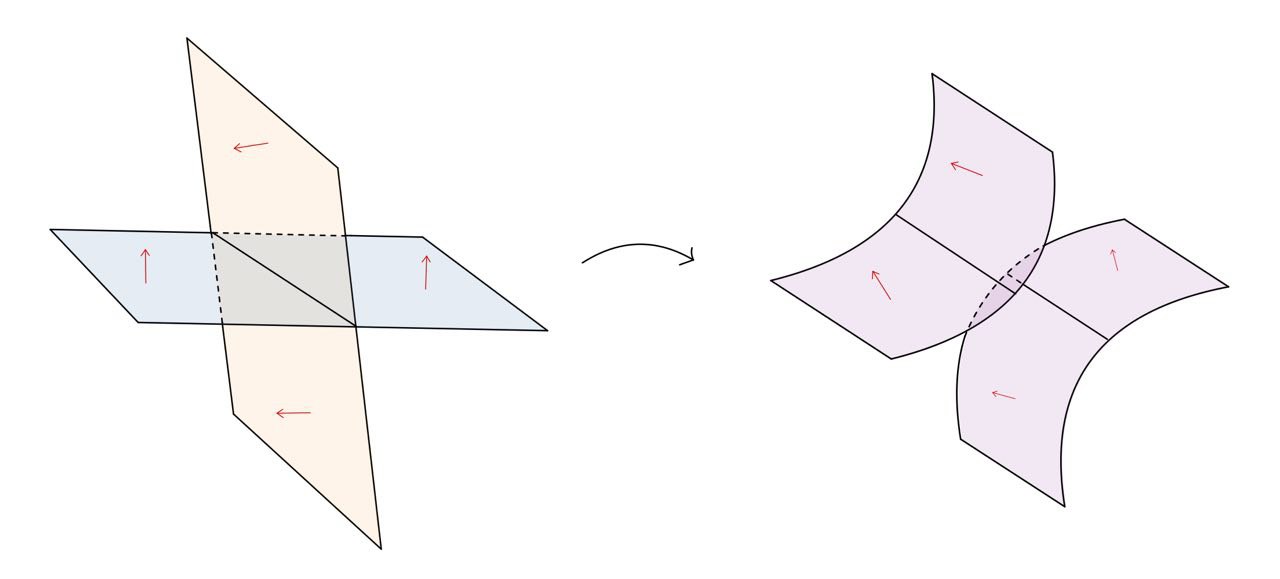}
    \caption{Local view of the oriented cut-and-paste operation of the surfaces $S_1$ and $S_2$. Red arrows indicate the normal orientation.}
    \label{fig: orientedcutandpaste}
\end{figure}

\bd(Oriented cut-and-paste) Let $S_1,S_2\subset M$ be properly embedded oriented surfaces, slightly isotoped to intersect transversely. There is only one oriented way to split $S_1$ and $S_2$ along $S_1\cap S_2$ and reglue back to obtain the oriented surface $S$, which locally looks as in Figure \ref{fig: orientedcutandpaste}. The surface $S$ is called the \emph{oriented cut-and-paste} (or the \emph{double-curve sum}) of $S_1$ and $S_2$. 
\ed

\bd(Thurston norm) Let $S\subset M$ be a properly embedded orientable surface. Denote with $\chi(S)$ the Euler characteristic of $S$. If $S$ is connected, call $$\chi_-(S):=\max (0,-\chi(S)).$$
If $S$ is the disjoint union of connected surfaces $S_1,...,S_k$, define $$\chi_-(S):=\chi_-(S_1)+...+\chi_-(S_k).$$
For every class $\alpha\in H_2(M,\partial M;\Z)$, let $$x(\alpha)=\min_{[S]=\alpha}\chi_-(S),$$ where $S\subset M$ varies among the properly embedded oriented surfaces representing $\alpha$.
\ed

\bt[\cite{norm}] The function $x:H_2(M,\partial M;\Z)\to \Z$ can be extended to a seminorm \\ $x:H_2(M,\partial M;\R)\to \R_{\ge 0}$, called the \emph{Thurston norm} associated to the manifold $M$.
\et
\noindent\emph{Sketch of the proof.} Thurston verifies that $x$ is $\N$-homogeneous on $H_2(M,\partial M;\Z)$ and that $\chi_-(S)\le \chi_-(S_1)+\chi_-(S_2)$ whenever $S$ is the oriented cut-and-paste of incompressible and $\partial$-incompressible surfaces $S_1$ and $S_2$. Then, $x$ can be extended by homogeneity to a convex function $x: H_2(M,\partial M;\Q)\to \Q$, which can further be continuously extended to a seminorm $x: H_2(M,\partial M;\R)\to \R_{\ge 0}$. 

\bigskip

To compute a (semi)norm on $\R^n$ it is enough to know its unit ball. A handy feature of Thurston norm is that its ball can be recovered by knowing a finite amount of data, thanks to the following result.

\bt[\cite{norm}]\label{thm: ball is a polyhedron} The unit ball $B_x$ of the Thurston norm $x$ on $M$ is a finite (but possibly unbounded) convex polyhedron, symmetric with respect to the origin, and whose vertices are rational points in $H_2(M,\partial M;\R)$. 
\et

When $x$ is not a norm, i.e. when some surface with non-negative Euler characteristic is not null-homologous, the null-space of $x$ is a vector space contained in $B_x$. Of course, $B_x$ is bounded if and only if $x$ is a norm on $H_2(M,\partial M;\R)$. 

When $M$ is a graph manifold, many of its second homology classes have null Thurston norm, as we will see. In this case, it is most useful to study a quotient of $H_2(M,\partial M;\R)$ where Thurston norm induces a norm.  

\bd[Nonvanishing Thurston norm and nonvanishing homology]\label{def: nonvanishing} Let $V\subset H_2(M,\partial M;\R)$ be the null-space of the Thurston seminorm $x$. We will refer to the quotient $$H_2^{nv}(M,\partial M;\R):=H_2(M,\partial M;\R)/V$$ as the \emph{nonvanishing real second homology group} of $M$ (with respect to the Thurston norm). Analogously, we can define the \emph{nonvanishing integral (resp. rational) second homology group} $H_2^{nv}(M,\partial M;\Z)$ (resp. $H_2^{nv}(M,\partial M;\Q)$) of $M$, as a quotient of $H_2(M,\partial M;\Z)$ (resp. $H_2(M,\partial M;\Q)$). We will think of $H_2^{nv}(M,\partial M;\Z)$ as the integral lattice in $H_2^{nv}(M,\partial M;\R)$.

The seminorm $x$ on $H_2(M,\partial M;\R)$ induces a well-defined norm $x^{nv}$ on $H_2^{nv}(M,\partial M;\R)$, by the formula $$x^{nv}(\alpha+V):=x(\alpha)$$ for every class $\alpha\in H_2(M,\partial M;\R)$. 

To simplify the notation, when not mentioned, we will understand the real coefficients for the homology: for example, we will refer to $H_2(M,\partial M;\R)$ just as $H_2(M,\partial M)$.
\ed

\br Let $B_x^{nv}$ be the unit ball in $H_2^{nv}(M,\partial M)$ with respect to the nonvanishing Thurston norm. As a consequence of Theorem \ref{thm: ball is a polyhedron}, $B_x^{nv}$ is a convex compact polyhedron symmetric through the origin. Moreover the polyhedron $B_x$ in $H_2(M,\partial M)$ has the structure of a product polyhedron $B_x^{nv}\times V$, where $V$ is the null-space of $x$.
\er

The combinatorial structure of the Thurston ball $B_x$ retains information about the ways $M$ fibers over the circle, if any.

\bt\label{thm: fibered faces}\cite{norm} If $S\subset M$ is the fiber of a fibration of $M$ over $S^1$, then $[S]\in H_2(M,\partial M)$ lies in the interior of the cone over a top-dimensional boundary face of the Thurston ball $B_x$. Viceversa, there is a collection of top-dimensional boundary faces of $B_x$ such that every integral class lying in the interior of the cone over one such face is represented by a fiber of a fibration of $M$ over $S^1$. 
\et

The open boundary faces of $B_x$ satisfying Theorem \ref{thm: fibered faces} are called \emph{fibered faces}. The cone over a fibered face is called a \emph{fibered cone} of the Thurston norm, and the classes lying in it are called \emph{fibered classes}.

\bigskip

An incompressible torus $T$ in $M$ allows us to focus on the manifold obtained by cutting $M$ along $T$. It turns out that the Thurston norm behaves well with respect to this decomposition.

\bl\label{lemma: additive norm}\cite[Proposition 11.2]{la21} Let $M$ be a compact oriented $3$-manifold with incompressible boundary. Let $T\subset M$ be a disjoint union of incompressible tori and let $M'$ be $M$ cut along $T$, seen as a submanifold of $M$. Consider a properly embedded oriented surface $S\subset M$, being transverse to $T$. Indicate with $x_M$ and $x_{M'}$ the Thurston norms on $H_2(M,\partial M)$ and $H_2(M',\partial M')$ respectively. The following equality holds $$x_M([S])=x_{M'}([S\cap M']).$$
\el

\subsection{Seifert fibered manifolds and their Thurston norm} 
We now introduce the building blocks for graph manifolds, which are Seifert fibered manifolds. For further details, see \cite{martelli} for example. 

For the whole section, let $\Sigma$ be a compact surface. Recall that, up to isomorphism, there is only one orientable circle-bundle over $\Sigma$, which we will indicate with $\Sigma\widetilde{\times}S^1$. 

\bd[Seifert fibered manifolds]\label{def: seifert manifolds} Fix an orientation on $\Sigma\widetilde{\times} S^1$. 

Given a circle fiber $c$ in $\Sigma\widetilde{\times} S^1$, let $T$ be the boundary of a small tubular neighbourhood $N(c)$ of this fiber. Slightly isotope the torus $T$ so that it is a union of circle fibers. We can choose a circle fiber as a longitude $l$ of $T$ and the intersection $T\cap(\Sigma\times{0})$ as a meridian $m$. Give orientations to $m$ and $l$ so that the couple $(m,l)$ is a positively oriented basis of $H_1(T;\Z)$, where $T$ has the boundary orientation of $\Sigma\widetilde{\times} S^1-N(c)$.

Consider some pairs of coprime integers $(p_1,q_1),...,(p_r,q_r)$ with $p_i\ne 0$ for every $i$. We call $(\Sigma, (p_1,q_1),...,$ $ (p_r,q_r))$ the manifold obtained from $\Sigma\widetilde{\times} S^1$ by Dehn surgery as follows. Choose $r$ disjoint circle fibers $c_1,...c_r$ and $r$ disjoint small circle-fibered tubular neighbourhoods $N(c_1),...,N(c_r)$ of them. Let $T_1,T_2,...,T_r$ be the boundary tori of these tubular neighbourhoods, each equipped with a meridian $m_i$ and a longitude $l_i$ as above. For every $i$, do a $(p_i,q_i)$-Dehn surgery on $T_i$: drill the interior of $N(c_i)$ out of $\Sigma\widetilde{\times} S^1$ and glue a solid torus back to $T_i$ so as to kill the slope $p_im_i+q_il_i$. 

The manifold $(\Sigma, (p_1,q_1),...,(p_r,q_r))$ can be endowed with a singular fibration by circles (see \cite{martelli}) over the base orbifold $(\Sigma,p_1,...,p_r)$, i.e. the orbifold whose base surface is $\Sigma$ with $r$ cone points, with orders $p_1,...,p_r$. This singular fibration is restricted to the regular circle bundle on $\Sigma\widetilde{\times} S^1-\left(\cup_i N(c_i)\right)$ and has singular fibers only corresponding to the cone points of $(\Sigma,p_1,...,p_r)$. This special type of singular fibration is called \emph{Seifert fibration}, hence we call $(\Sigma, (p_1,q_1),...,(p_r,q_r))$ a \emph{Seifert fibered manifold} (over the base orbifold $(\Sigma,p_1,...,p_r)$).
\ed

The classification of Seifert fibered manifolds up to homeomorphism is completely understood \cite{martelli}. A weaker problem concerns the classification of Seifert fibered manifolds up to \emph{isomorphism}. If we are given two Seifert fibered manifolds $M$ and $M'$ with fixed Seifert fibrations on them, we say they are isomorphic if there is a diffeomorphism between $M$ and $M'$ respecting the fibrations. The same manifold can admit two different (non-isomorphic) Seifert fibrations. 

However, in most cases there will be no ambiguity, e.g. for manifolds with boundary, thanks to the following proposition.

\bprop\cite[Theorem 10.4.19]{martelli} If $\partial\Sigma\ne\emptyset$ and $\Sigma$ is not a disc nor a Moebius band, then $(\Sigma, (p_1,q_1),...,(p_r,q_r))$ admits only one Seifert fibration up to isomorphism.
\eprop

To avoid ambiguity for the cases not covered by the previous proposition, we will understand that the manifold $(\Sigma, (p_1,q_1),...,(p_r,q_r))$ is endowed with the Seifert fibration coming by the construction of Definition \ref{def: seifert manifolds}.

\bigskip

When dealing with Seifert fibered manifolds, two important invariants are always concerned, namely the base orbifold Euler characteristic and the Euler number. Here, we will use the conventions of \cite{neumann}.

\bd[Euler invariants]\label{def: invariants} Let $M=(\Sigma, (p_1,q_1),...,(p_r,q_r))$ be a Seifert fibered manifold. The \emph{base orbifold Euler characteristic} of $M$ is defined as $$\chi^{\text{orb}}(M):=\chi(\Sigma)-\sum_{i=1}^r\left(1- \frac 1{|p_i|}\right).$$
When $M$ is closed (i.e. when $\Sigma$ is), the \emph{Euler number} of $M$ is defined as $$e(M):=-\sum_{i=1}^r \frac {q_i}{p_i}.$$
When $M$ has nonempty boundary, we reduce to the closed case by Dehn-filling along the boundary components as follows. A \emph{system of meridians} for $M$ is a choice of a simple closed curve non-isotopic to a fiber on each boundary torus of $M$. Fix a system of meridians for $M$. Let $\overline M$ be the closed manifold obtained from $M$ by Dehn-filling along its boundary tori so as to kill each slope in the system of meridians. The Seifert fibration on $M$ extends to a Seifert fibration on $\overline M$ and we call $e(M):=e(\overline M)$ the \emph{Euler number} of $M$ with respect to the chosen system of meridians.
\ed

\noindent We can think of $\chi^{\text{orb}}(M)$ as a measure of the topological complexity $\chi^{\text{orb}}(\Sigma)$ of the base orbifold. The defining formula is suitably chosen to be multiplicative under finite orbifold coverings, which is also an ingredient used in the proof of Proposition \ref{prop: norm for seifert}.

For a closed Seifert fibered manifold $M$, the Euler number is zero exactly when the Seifert fibration is finitely covered by a trivial circle bundle (see \cite[Corollary 10.3.29]{martelli}). We can then consider the Euler number as a measure of the twistedness of the circle fibers in the Seifert fibration.

\bigskip

For the remaining part of this section, let $M=(\Sigma, (p_1,q_1),...,(p_r,q_r))$ with the induced Seifert fibration. The fibration in $M$ helps us to describe the surfaces sitting inside $M$. 

\bd[Vertical and horizontal surfaces] Let $S\subset M$ be a properly embedded surface. We say that $S$ is \emph{vertical} if it is a union of regular circle fibers. We say that $S$ is \emph{horizontal} if it intersects the regular fibers transversely.
\ed

\bprop\cite[Propositions 10.4.9 and 10.4.10]{martelli} Let $S\subset M$ be a properly embedded surface and suppose $M$ is not diffeomorphic to $S^2\times S^1$ nor to $\mathbb{RP}^2\widetilde{\times} S^1$. If $S$ is incompressible, $\partial$-incompressible and no component of $S$ is $\partial$-parallel, then $S$ is isotopic to either a vertical or a horizontal surface. 
\eprop

We can now describe the nonvanishing second homology of a Seifert fibered manifold.

\bprop\label{prop: norm for seifert} Suppose that $M$ is connected. The nonvanishing second homology of $M$ is generated by horizontal surfaces and $$H_2^{nv}(M,\partial M)\cong \begin{cases}
    \R\text{ if $\chi^{orb}(M)<0$ and $\partial M\ne\emptyset $}\\ \R\text{ if $\chi^{orb}(M)<0$ and $e(M)=0$}\\ 0\text{ otherwise.}
\end{cases}$$
Additionally, given any regular fiber $\ell$ of the Seifert fibration, the following formula holds for every $\alpha\in H_2(M,\partial M)$: 
\begin{align}\label{equation: norm of Seifert} x(\alpha)=\chi_-^{orb}(M)\cdot |i(\alpha, \ell)|=\chi_-^{orb}(M)\cdot|\psi(\alpha)|,
\end{align}
where $i(\alpha,\ell)$ denotes the algebraic intersection number, $\psi\in H^2(M)\cong\hom(H_2(M),\R)$ is the Poincare dual of $\ell$ and $\chi_-^{orb}(M)=\max\left(0, -\chi^{orb}(M)\right).$
\eprop
\bp Let $S\subset M$ be a properly embedded oriented surface. Up to homology, we can assume $S$ is incompressible and $\partial$-incompressible. In particular, after possibly discarding null-homologous components, we can suppose $S$ is either vertical or horizontal. 

If $S$ is vertical then it must be either an annulus or a torus, so $\chi(S)=0$ and $[S]$ belongs to the null-space of the Thurston norm.

As proved in \cite[Proposition 10.4.8]{martelli}, if $M$ is closed with $e(M)\ne 0$, then $M$ contains no horizontal surface, so $H_2^{nv}(M,\partial M)=0$. So, assume that $\partial M\ne \emptyset$ or $e(M)=0$, in which case $M$ admits a horizontal surface.

Suppose that $S$ is horizontal and connected. Then, the projection $M\to \Sigma$ restricts to a degree-$d$ orbifold covering $S\to\Sigma$, where $d$ is the number of times that $S$ transversely intersects each fiber. Thus, $$\chi(S)=\chi^{orb}(S)=d\chi^{orb}(\Sigma)=d\chi^{orb}(M).$$
Observe that since $S$ is connected, it must intersect every regular fiber always with the same sign. Indeed, the set of points intersecting the regular fibers positively is open in $S$, as well as the set of points intersecting them negatively, and these two sets together cover $S$.
In particular, given any regular fiber $\ell$ of the Seifert fibration, \begin{equation}\label{equation: degree}
    d=|S\cap\ell|=|i([S],\ell)|.
\end{equation}
By cutting $M$ along $S$, we obtain an oriented interval bundle over $S$, hence $S\times [0,1]$. So, the manifold $M$ fibers over the circle with fiber $S$, and this implies that $S$ is norm-minimising in $H_2(M,\partial M)$ \cite{norm}, i.e. \begin{align}\label{equation: proof of norm of seifert}x([S])=\chi_-(S)=d\max (0,-\chi^{orb}(M)).\end{align}

Finally, suppose $S'$ is another oriented horizontal surface, transversely intersecting each fiber $d'$ times. After possibly reversing the orientation of $S'$, we can suppose that the surface given by the oriented cut-and-paste $\overline S$ of $d'$ parallel copies of $S$ and $d$ parallel copies of $S'$ has algebraic intersection $0$ with $\ell$. If $\overline S$ is homologous to a horizontal surface, such a horizontal surface consists of parallel connected components (since each component is the fiber of a fibration of $M$ over the circle), and thus $i([\overline S],\ell)=0$ implies that $\overline S$ is null-homologous. If $\overline S$ is homologous to a vertical surface, then it annihilates the Thurston norm. This proves that the nonvanishing homology has dimension at most one, the reasoning above shows exactly when it vanishes and formulae (\ref{equation: degree}) and (\ref{equation: proof of norm of seifert}) together imply formula (\ref{equation: norm of Seifert}). 
\ep

Once we know the nonvanishing homology of $M$, we just need to understand the null-space of the Thurston norm to have the whole picture of the Thurston norm on $H_2(M,\partial M).$

\bprop\label{prop: vertical homology} Suppose $\Sigma$ is orientable and connected. The subspace of $H_2(M,\partial M)$ generated by properly embedded oriented vertical surfaces is isomorphic to $H_1(\Sigma,\partial \Sigma)$ (the singular homology of the underlying topological surface $\Sigma$). The subspace of $H_2(M,\partial M)$ generated by closed oriented vertical surfaces is isomorphic to the image of $H_1(\Sigma)$ in $H_1(\Sigma,\partial \Sigma)$ and has dimension equal to twice the genus of $\Sigma$.
\eprop
\bp Consider the projection map $\pi:M\to \Sigma$ on the base orbifold. The induced map $\pi^*:H^1(\Sigma)\to H^1(M)$ through Lefschetz Duality gives $\pi^*:H_1(\Sigma,\partial \Sigma)\to H_2(M,\partial M).$ If $\mu\subset \Sigma$ is a properly embedded multicurve, slightly moved so as not to touch the cone points of $\Sigma$, it is easy to see that $\pi^*([\mu])=[\pi^{-1}(\mu)= \mu\times S^1]$. Since the image of $\pi^*$ is the subspace generated by vertical surfaces, we want to show that $\pi^*$ is injective.

The manifold $M$ is obtained by non-fiber-parallel Dehn-filling on some product manifold $M'\cong\Sigma'\times S^1$. Here, $\Sigma$ comes from capping off some boundary component of $\Sigma'$ with discs.  Consider the commutative diagrams:

\begin{center}
    % https://tikzcd.yichuanshen.de/#N4Igdg9gJgpgziAXAbVABwnAlgFyxMJZABgBpiBdUkANwEMAbAVxiRAFkByEAX1PUy58hFAEZyVWoxZt2vfiAzY8BImVGT6zVohAAdPQGUsAcwC2dbnwHLhRcRupaZug8fN15NoapQAmCSdpHRAACQA9UQAKLgBKL0VBFRFkAGZAqW02COj2eOtE219kAMdMlzDIqLdTC058hSUfFPSy5xCc6qNaunzJGCgTeCJQADMAJwgzJDIQHAgkAPKQgzQsKwUJqaRxOYXEdLm6LAY2SDBWIKzXPTWErenEWfmd6hxj090ACwgIAGsQFcKgAre6TR5LF4HN4fNg-f6A5ZsLBg7aIABsb32ABYgSEumt6uEAFSIhh0ABGMAYAAUknZdONTF8cKjHgB2LFIACseLYqywJLJlOpdKKIhATJMLLZPK5iFxSN0wKF1HJVNp9N8kuZrIKDyQnL2SExSpAgtJapFmvFbClMp4FB4QA
\begin{tikzcd}
M' \arrow[d, "\pi'"] \arrow[r, "j", hook] & M \arrow[d, "\pi"] & H^1(M')                             & H^1(M) \arrow[l, "j^*"']                          \\
\Sigma' \arrow[r, "i", hook]              & \Sigma             & H^1(\Sigma') \arrow[u, "(\pi')^*"'] & H^1(\Sigma) \arrow[u, "\pi^*"'] \arrow[l, "i^*"']
\end{tikzcd}
\end{center}
the maps $i,j$ indicating the inclusions, $\pi'$ the projection induced by the product fibration.
The map $(\pi')^*$ is injective with left inverse the map $H^1(M')\to H^1(\Sigma')$ induced by the inclusion of a leaf $\Sigma'\to M'$. The map $i^*$ is injective too, because $\Sigma$ is obtained from $\Sigma'$ by just capping off with discs. By the commutativity of the diagram on the right, we conclude that $\pi^*$ is injective.

Now, it easily follows that the subspace of $H_2(M,\partial M)$ generated by closed oriented vertical surfaces is isomorphic to the image of $H_1(\Sigma)$ in $H_1(\Sigma,\partial \Sigma)$. The only remaining fact to check is the one about the rank of this subspace. If $\partial \Sigma=\emptyset$ there is nothing to show. Suppose then that $\Sigma$ has $b\ge1$ boundary components and genus $g$. Consider the long exact sequence of the couple $(\Sigma,\partial \Sigma)$ for reduced homology: $$H_1(\Sigma)\to H_1(\Sigma,\partial \Sigma)\to \widetilde H_0(\partial\Sigma)\to \widetilde H_0(\Sigma).$$ By Lefschetz Duality, $H_1(\Sigma,\partial\Sigma)\cong H_1(\Sigma)\cong \R^{2g+b-1}$, so the exact sequence becomes $$H_1(\Sigma)\to \R^{2g+b-1}\to \R^{b-1}\to 0.$$ In particular, $$\dim \im\left(H_1(\Sigma)\to H_1(\Sigma,\partial \Sigma)\right)=\dim\Ker \left(\R^{2g+b-1}\to \R^{b-1}\to 0\right)=2g+b-1-(b-1)=2g.$$
\ep

\subsection{Graph manifolds and their reduced plumbing matrices}\label{subsec: graph manifolds}
In this section, we introduce graph manifolds and some results by Neumann \cite{neumann} that will be essential for our ends. We will be interested only in the closed orientable case, even if the stated facts have broader generality. 

\bd[Graph manifold] Let $M$ be a $3$-dimensional closed oriented manifold. We say that $M$ is a \emph{graph manifold}, if there is a collection of tori $T_1,...,T_k$ in $M$ such that $M$ cut along this collection is a disjoint union of Seifert fibered manifolds. If the Seifert fibered pieces all have orientable base orbifolds, we say that $M$ is a \emph{good} graph manifold. Notice that if $M$ is a non-good graph manifold, then $M$ is double-covered by a good graph manifold.

When $M$ is a graph manifold, there is a canonical family of decomposing tori in $M$ up to isotopy. This family is characterized as the minimal one (with respect to the inclusion) cutting $M$ into Seifert fibered pieces, and the decomposition is the \emph{JSJ decomposition} of $M$ \cite{martelli}. When working with a graph manifold $M$ and unless otherwise mentioned, we will always consider its subdivision into Seifert fibered pieces coming from the JSJ decomposition.

Suppose that the JSJ decomposition of $M$ consists of the Seifert fibered pieces $M_1,...,M_n$. We can build $M$ back from $M_1,...,M_n$ by means of the \emph{decomposition graph} $\Gamma$ of $M$. The vertices of $\Gamma$ are in $1$-$1$ correspondence with the pieces $M_1,...,M_n$, and there is an edge connecting the vertices $i$ and $j$ for every decomposing torus separating $M_i$ from $M_j$. Then, each edge is endowed with a self-diffeomorphism of the torus, collecting information on how to glue the two pieces back along their common boundary.

Observe that in our notation, a disconnected graph represents a disconnected graph manifold. 

\bigskip

To describe the Thurston norm of $M$, we will need less information than the whole data above. For every $i$, fix an orientation on each piece $M_i$, as well as a Seifert fibration on it. Call $\chi_i:=\chi^{orb}(M_i)$ and $e_i:=e(M_i)$. Here, the Euler number is computed with respect to the system of meridians given by choosing as a meridian on each boundary torus $T$ of $M_i$ a regular fiber coming from the other piece $M_j$ across $T$. For every edge $E$ connecting the vertices $i$ and $j$ in $\Gamma$ (we will write $iEj$), call $p(E)$ the algebraic intersection number between the two regular fibers of $M_i$ and $M_j$ on the separating torus represented by $E$.

We will call \emph{simplified decomposition graph} of $M$ the graph with the same vertices and edges as $\Gamma$, where we label each vertex $i$ with the couple $(e_i,\chi_i)$ and each edge $E$ with the integer $p(E)$.

Finally, the \emph{reduced plumbing matrix} of $M$ is the matrix $A=(a_{i,j})_{i,j=1,...,n}$ with $$a_{i,j}:=\begin{cases}
    e_i+2\sum_{iEi}\frac 1{p(E)} \text{ if $i=j$} \\
    \sum_{iEj}\frac 1{p(E)} \text{ if $i\ne j$}.
\end{cases}$$
\ed

Observe that the number $p(E)$ is nonzero because the pieces $M_1,...,M_n$ come from the JSJ decomposition of $M$.

As discovered by Neumann, the reduced plumbing matrix yields much information about how the surfaces sit inside the ambient graph manifold. Notice that the following is a slightly corrected version of the cited result by Neumann.

\bl\cite[Lemma 4.2]{neumann}\label{lemma: arising tuples} Let $M$ be closed oriented good graph manifold with pieces $M_1,...,M_n$. Assume that $M$ is not a torus-bundle over $S^1$. Let $(l_1,...,l_n)$ be a tuple of integers. The following properties are equivalent:
\begin{itemize}
    \item There are a natural $n_0>0$ and a properly embedded oriented surface $S\subset M$ intersecting the regular fibers of $M_i$ algebraically $n_0l_i$ times, for each $i=1,...,n$;
    \item The vector $(l_1,...,l_n)$ annihilates the reduced plumbing matrix of $M$.
\end{itemize}
\el 

\bigskip

\noindent In the original version of \cite[Lemma 4.2]{neumann}, in the place of the first property there was \begin{itemize}
    \item There is a properly embedded oriented surface $S\subset M$ intersecting the regular fibers of $M_i$ algebraically $l_i$ times, for each $i=1,...,n$.
\end{itemize}
This stronger statement is false in general: as an example, consider a closed Seifert fibered manifold $N$ with Euler number $0$ and nonintegral orbifold Euler characteristic $q/p$. A horizontal surface of $N$ covers the base orbifold, the degree of the cover being the algebraic intersection $\ell$ between the fiber and the horizontal surface. Since the horizontal surface has integral Euler characteristic $|\ell| q/p$, the integer $p$ necessarily divides $\ell$. Hence, even though the reduced plumbing matrix for $N$ is $0:\Z\to\Z$, the set of the realizable values for $\ell$ is contained in $p\Z$.

Observe that the corrected version of Lemma \ref{lemma: arising tuples} is enough to prove the results of \cite{neumann}, by just replacing the original lemma. 

To justify the correction, we now give an idea of Neumann's proof, which still proves the corrected version of Lemma \ref{lemma: arising tuples}. 

\bigskip

\noindent\textit{Sketch of the proof.} In \cite{neumann}, Neumann considers $M$ as a \emph{plumbing manifold}, often relying on the results of \cite{plumbing}. A plumbing manifold is a $3$-manifold $N$ that admits a family of decomposing tori $\mathcal T$ such that:
\begin{itemize}
    \item[(i)] When cut along $\mathcal T$, $N$ gives a disjoint union of circle-bundle pieces $N_1,...,N_h$;
    \item[(ii)] If two pieces $N_i$ and $N_j$ are adjacent in $N$ to some torus $T\in\mathcal T$, then the gluing map has the form $\pm\begin{pmatrix}
    0 & 1 \\ 1&0
\end{pmatrix}$, with respect to a choice of meridians and longitudes induced by the fibrations.
\end{itemize}
A Seifert fibered manifold admits such a decomposition, and so does a graph manifold. In particular, a plumbing decomposition of $M$ is a refinement of its JSJ decomposition. Equivalently, the JSJ tori of $M$ are a subset of any family of plumbing tori $\mathcal T$ for $M$. 

Given a plumbing decomposition of $M$, we can define a plumbing graph $\Gamma'$, which is constructed analogously to the decomposition graph of a graph manifold. Every vertex of $\Gamma'$ corresponds with a plumbing piece of $M$, and every edge indicates that the relative two pieces are glued along a torus boundary. The graph $\Gamma'$ is decorated with labels on vertices indicating the topology of the pieces' base surfaces (genus, number of boundary components) and their Euler numbers. A sign $+$ or $-$ on every edge indicates which of the two possible glueing maps stands in the relative glueing. 

A useful feature of the plumbing graph $\Gamma'$ is that it can be used to build up a compact orientable $4$-manifold $X$ whose boundary is $M$. Essentially, vertices of $\Gamma'$ codify for disc-bundles, and edges of $\Gamma'$ tell how to glue $2$-handles to these bundles. 

Neumann then considers the following part of the cohomological long exact sequence of the couple $(X,M)$: $$H^1(M;\Z)\to H^2(X,M;\Z)\to H^2(X;\Z).$$ By Lefschetz Duality and the Universal Coefficients Theorem, it can be rewritten as $$H_2(M;\Z)\overset{\alpha}{\longrightarrow} H_2(X;\Z)\overset{\beta}{\longrightarrow} \hom (H_2(X;\Z),\Z).$$
The base surfaces of the plumbing pieces $M_1',...,M_k'$ of $M$ give a basis of $H_2(X;\Z)\cong \Z^k$. In such basis, $\alpha$ sends the class of a properly embedded oriented surface $S\subset M$ to the tuple $(\ell_1',...,\ell_k')$, where $\ell_i'$ is the algebraic intersection of $S$ with the fiber of $M_i'$. The map $\beta$ is the intersection form sending a surface $S\subset X$ to the algebraic intersection homomorphism $\langle [S],\cdot\rangle$, represented by a matrix $\mathcal S$ in the aforementioned basis. Now, we just need to establish the relation between the tuples $(\ell_1',...,\ell_k')$ and $(\ell_1,...,\ell_n)$, as well as the relation between $\mathcal S$ and the reducing plumbing matrix $A$ of $M$. 

Each piece $M_j$ of the JSJ decomposition of $M$ corresponds to a subgraph $\Gamma_j$ of $\Gamma'$, and the subgraphs obtained in this way all together cover $\Gamma'$, while pairwise intersecting just in a union of vertices. In every subgraph $\Gamma_j$, there is a preferred vertex, called the \emph{node} of the subgraph. The node is a vertex with positive genus or with degree $\ge 3$.  Label the vertices of $\Gamma'$ so that the nodes are labelled $1,...,n$. The matrix $\mathcal S$ can be modified through a sequence of column and row moves to get the direct sum of $A$ with a diagonal $(k-n)\times(k-n)$ matrix. Each row (resp. column) move consists of adding to a row (resp. column) a rational multiple of one of the last $k-n$ rows (resp. columns). In particular, if $(\ell_1',...,\ell_k')\in\Z^k$ annihilates $\mathcal S$, then $(\ell_1',...,\ell_n')$ annihilates $A$. Viceversa, if $(\ell_1,...,\ell_n)\in\Z^n$ annihilates $A$, then there are $\ell_{n+1}',...,\ell_k'\in\Q$ such that $(\ell_1,...,\ell_n,\ell_{n+1}',...,\ell_k')$ annihilates $\mathcal S$. Since $\ell_{n+1}',...,\ell_k'$ can be nonintegral, we could need to multiply $(\ell_1,...,\ell_n,\ell_{n+1}',...,\ell_k')$ by a nonzero natural number $n_0$ to obtain a tuple of integers annihilating $\mathcal S$.

\section{The Thurston norm of a graph manifold}\label{sec: norm of a graph manifold}
This section aims to compute the Thurston norm of a closed oriented graph manifold $M$. 

Let $M_1,...,M_n$ be the pieces of the JSJ decomposition of $M$, $\chi_i$ and $e_i$ be their associated invariants, and $A$ be the reduced plumbing matrix of $M$, seen as a linear map $\R^n\to\R^n$. In order to have nontrivial nonvanishing norm on $M_i$, we will assume that $\chi_i<0$ for every $i=1,...,n$.

For every $i=1,...,n$, let $\psi_i\in H^2(M)\cong\hom(H_2(M),\R)$ be the Poincare dual of a regular fiber $l_i\subset M_i$. Finally, define the linear map $\psi=(\psi_1,...,\psi_n):H_2(M)\to\R^n$.

\bl\label{lemma: norm} For every $\alpha\in H_2(M)$, the following formula holds: \begin{equation}\label{equation: norm}
x(\alpha)=-\sum_{i=1}^n \chi_i |\psi_i(\alpha)|.\end{equation}
\el
\bp It is sufficient to prove the identity (\ref{equation: norm}) for integral classes. Thanks to Lemma \ref{lemma: additive norm}, for every properly embedded oriented surface $S\subset M$, after possibly a slight isotopy to make $S$ transverse with the decomposing tori of $M$, we have $$x([S])=\sum_{i=1}^n x_{M_i}([S\cap M_i]),$$ where $x_{M_i}$ is the Thurston norm on $H_2(M_i,\partial M_i)$. Each $M_i$ is a Seifert fibered manifold and equality (\ref{equation: norm of Seifert}) concludes the proof.
\ep

\bprop\label{prop: computation} If $M$ is good, then the following sequence of vector spaces is exact $$0\rightarrow V\rightarrow H_2(M)\overset{\psi}{\longrightarrow}\R^n\overset{A}{\longrightarrow}\R^n$$ where $V$ is the null-space of the Thurston norm on $H_2(M)$. In particular, the nonvanishing homology $H_2^{nv}(M)$ is isomorphic to $\ker A\subset \R^n$ through an isomorphism sending $H_2^{nv}(M;\Q)$ onto $\ker A\cap \Q^n$. Hence, a rational basis of  $H_2^{nv}(M)$ corresponds through $\psi$ with a rational basis of $\ker A$.
\eprop
\bp It is enough to show that the sequence of $\Q$-vector spaces $$0\rightarrow V\cap H_2(M;\Q) \rightarrow H_2(M;\Q)\overset{\psi}{\longrightarrow}\Q^n\overset{A}{\longrightarrow}\Q^n$$ is exact. Since $\chi_i<0$ for every $i$, $M$ cannot be a torus-bundle. Indeed, a non-Seifert fibered torus-bundle has JSJ decomposition equal to a product (torus)$\times [0,1]=$(annulus)$\times S^1$, given by cutting along a torus fiber. Hence, Lemma \ref{lemma: arising tuples} implies that $\psi(H_2(M;\Q))=\ker A\cap\Q^n$. Finally, given that $\chi_i<0$ for each $i$, equality (\ref{equation: norm}) implies $V\cap H_2(M,\Q)=\ker \psi$.
\ep

\bprop\label{prop: nonvanishing norm} Suppose $M$ is good. Let $v_1,...,v_d\in \R^n$ be a basis of $\ker A$. Consider the $n\times d$ matrix $P$ with columns respectively $v_1,...,v_d$ and let $\beta_1,...,\beta_n\in\R^d$ be the rows of $P$. In the coordinates with respect to $(v_1,...,v_d)$, the nonvanishing Thurston norm on $\ker A$ respects the formula: \begin{equation}\label{equation: nonvanishing norm}
    x^{nv}(v)=-\sum_{i=1}^n\chi_i |\langle \beta_i,v\rangle|
\end{equation}  for every $v\in\R^d$. 
\eprop
\bp  Let $\alpha_1,...,\alpha_d$ be the basis of $H_2^{nv}(M)$ corresponding with $v_1,...,v_d$ through $\psi$, i.e. such that $\psi(\alpha_j)=v_j$ for every $j$. 

For any $i=1,...,n$ and $a_1,...,a_d\in\R$, $$\psi_i(a_1\alpha_1+...+a_d\alpha_d)=a_1\psi_i(\alpha_1)+...+a_d\psi_i(\alpha_d)=a_1(P)_{i1}+...+a_d(P)_{id}=\langle \beta_i, (a_1,...,a_d) \rangle.$$
Now, formula (\ref{equation: norm}) implies the thesis.
\ep

Proposition \ref{prop: nonvanishing norm} expresses $x^{nv}$ as a sum of absolute values of linear functionals. Construction \ref{construction} describes the unit polyhedron of such a norm.

\bigskip

Let us investigate the null-space of the Thurston norm now.

\bprop\label{prop: dimension} Suppose $M$ is good and let $\Gamma$ be the decomposition graph of $M$. For $i=1,...,n$, let $g_i$ be the genus of the base orbifold of the piece $M_i$. The null-space of the Thurston norm on $H_2(M)$ has dimension $b_1(\Gamma)+2\sum_{i=1}^n g_i=|\Gamma|-\chi(\Gamma)+2\sum_{i=1}^n g_i$, where $|\Gamma|$ is the number of connected components of $\Gamma$ and $\chi(\Gamma)$ is its Euler characteristic. The second Betti number of $M$ is \begin{align*}
    b_2(M)=\dim\ker A+b_1(\Gamma)+2\sum_{i=1}^n g_i.
\end{align*}
\eprop
\bp Thanks to Lemma \ref{lemma: norm}, if $S\subset M$ is an oriented and properly embedded surface of norm zero, then $S\cap M_i$ is isotopic to a vertical surface in $M_i$ for every $i$. Since the regular fibers of two adjacent pieces of $M$ are never parallel, $S$ is a disjoint union of closed vertical surfaces in some $M_i$s. For each $i$, let $\overline{V_i}\subset H_2(M_i,\partial M_i)$ be the subspace generated by the classes of vertical closed surfaces in $M_i$. Choose a subspace $V_i\subset H_2(M_i)$ such that the map $H_2(M_i)\to H_2(M_i,\partial M_i)$ restricts to an isomorphism between $V_i$ and $\overline{V_i}$. Identify $V_i$ with its image under the inclusion $H_2(M_i)\subset H_2(M)$. Let $\mathcal T\subset H_2(M)$ be the subspace generated by the classes of JSJ tori of $M$. The null-space of the Thurston norm splits as the direct sum $\mathcal T\oplus V_1\oplus...\oplus V_n$. If we show that the dimension of $\mathcal T$ is $b_1(\Gamma)$, then the statement follows from Propositions \ref{prop: vertical homology} and \ref{prop: computation}. 

Let $E(\Gamma)$ and $V(\Gamma)$ be the sets of edges and vertices of $\Gamma$ respectively. Let $G_E$ and $G_V$ be the real vector spaces with bases $E(\Gamma)$ and $V(\Gamma)$ respectively. Choose an arbitrary orientation on each edge of $\Gamma$, so that each edge has now a \emph{source} and a \emph{sink}. For every $e\in E(\Gamma)$ and $v\in V(\Gamma)$, define the \emph{index} of $e$ at $v$ as $i(e,v)=0$ if none or both of the vertices of $e$ are $v$, otherwise $i(e,v)=+1$ if $v$ is the sink of $e$, and $i(e,v)=-1$ if $v$ is the source of $e$. Consider the linear map $\phi:G_V\to G_E$ such that, for every $v\in V(\Gamma)$, $$\phi(v)= \sum_{e\in E(\Gamma)}i(e,v)e.$$ There is a natural epimorphism $p:G_E\to \mathcal T$, associating to each edge $e$ the respective decomposing torus, with normal orientation pointing outward of the piece of $M$ corresponding to the source of $e$. Observe that we have an exact sequence $$G_V\overset{\phi}{\longrightarrow}G_E\overset{p}{\longrightarrow}\mathcal T\rightarrow 0.$$ The fact that $\im\phi\subset \Ker p$ follows from the fact that $\partial M_i$ consists of the tori corresponding to the edges $e$ at the vertex $v_i$ of $M_i$, with coorientation given by $i(e,v_i)$. On the other hand, if $W$ is a cobordism between some JSJ tori of $M$, then $W$ is a union of some $M_i$s. So, the $3$-chain represented by $W$ is a sum of the $3$-chains represented by some $M_i$s, hence $\im\phi\supset \Ker p$.

Finally, observe that the CW-chain complex associated with $\Gamma$ is given by the sequence $$0\rightarrow G_E\overset{\phi^t}{\longrightarrow}G_V\rightarrow 0$$ where the transpose $\phi^t: G_E\to G_V$ is taken with respect to the bases $E(\Gamma)$ and $V(\Gamma)$. Thus $$\dim \mathcal T= \dim G_E-\rank \phi=\dim G_E-\rank\phi^t=\dim \Ker\phi^t=\dim H_1(\Gamma)=b_1(\Gamma).$$
\ep

An interesting feature of the Thurston norm of graph manifolds is that either none or every top-dimensional face of the Thurston ball is fibered.

\bprop\label{prop: fibering} If $M$ fibers over the circle, then each top-dimensional face of the Thurston norm is fibered.
\eprop
\bp We first show that a class $\alpha\in H_2(M;\Z)$ is fibered if and only if $\psi_i(\alpha)\ne 0$ for each $i=1,...,n$. Suppose that $\psi_i(\alpha)\ne 0$ for each $i$. Without loss of generality, assume that $M$ is connected and $\alpha$ is primitive, so that $\alpha$ can be represented by a connected properly embedded oriented surface $S\subset M$. We can also assume that $S\cap M_i$ is a (possibly disconnected) horizontal surface in $M_i$. Cut $M$ open along $S$, then call $N$ the new manifold. Notice that the manifold $N_i$ obtained by cutting $M_i$ open along $S\cap M_i$ is an oriented interval-bundle over $S\cap M_i$, hence necessarily $(S\cap M_i)\times [0,1]$. Properly speaking, $N_i$ is a \emph{manifold with corners}, $\partial N_i$ being the union of the ``vertical boundary" $(S\cap \partial M_i)\times [0,1]$ and the ``horizontal boundary" $(S\cap M_i)\times \{0,1\}$. Since $\alpha\ne 0$, $S$ is non-separating and $N$ is the connected manifold obtained by gluing the interval bundles $(S\cap M_i)\times [0,1]$ along their ``vertical" boundary components. Each vertical boundary component is an interval-foliated annulus. Since the foliations of the annulus by intervals are all isotopic with each other, we can isotope the interval-bundles in a collar neighborhood of the vertical boundary of every $N_i$ so that they patch together to an interval-bundle of $N$ over the surface $S$. Again, this interval-bundle must be $S\times [0,1]$, thus $M$ fibers over the circle with fiber $S$.

For the vice-versa, suppose that $M$ fibers over the circle with fiber $S$. As in the proof of Proposition \ref{prop: computation}, $S$ cannot be a torus. In particular, since the fibration is a taut foliation, by \cite{norm} each JSJ torus of $M$ can be slightly isotoped to be transverse to the fibers of the fibration, apart for a finite number of saddle singularities. In fact, by the Poincare-Hopf Theorem, each JSJ torus can be taken to be entirely transverse to the fibration. Hence, the fibration restricts to a fibration onto every piece $M_i$. If any $S\cap M_i$ were vertical, then $M_i$ would be a torus-bundle or an annulus-bundle, which is excluded since $\chi_i<0$. It follows that $\psi_i([S])\ne 0$, for every $i=1,...,n$.

\bigskip

Having shown the claim, the proof of the proposition easily follows. Indeed, since $M$ fibers over the circle, there is a fibering class $\alpha_0\in H_2(M; \Z)$. Given $\alpha\in H_2(M;\Z)$ there are just finitely many $n\in \N$ such that, for some $i$, $\psi_i(n\alpha+\alpha_0)=0$. Thus, for $n$ big enough, the class $\alpha+\frac 1n\alpha_0$ is fibered. This shows that the fibered classes are dense in $H_2(M;\Z)$. Since a top-dimensional cone over the Thurston ball is open, containing either just fibered classes or no fibered classes at all, then all the top-dimensional cones over the Thurston ball must be fibered.
\ep

\section{Characterization of the norms that can be realised}\label{sec: characterization}
The formula (\ref{equation: norm}) implies that the nonvanishing Thurston norm of a closed oriented graph manifold can be expressed as a sum of absolute values of linear functionals $H_2^{nv}(M)\to \R$. Notice that those linear functionals are rational on $H_2^{nv}(M;\Q)$. In this section, we will prove that every norm with these properties can be realised as the nonvanishing Thurston norm on some closed oriented good graph manifold.

\bt\label{thm: every sum is graph} Let $||\cdot||$ be a norm on $\R^d$ that can be written as $$||\cdot ||:=\sum_{i=1}^n |\langle \beta_i,\cdot\rangle|$$ for some $\beta_1,...,\beta_n\in \Q^d$, and let $g_1,...,g_n$ be nonnegative integers. There exist a closed oriented good graph manifold $M$ with base surfaces of genera $g_1,...,g_n$ and a rational basis on $H_2(M)$, such that the nonvanishing Thurston norm on $H_2^{nv}(M)$ in the induced basis coincides with $||\cdot||$. Moreover, we can ask for either all or none of the top-dimensional Thurston cones of $H_2(M)$ to be fibered.
\et

\bp First observe that the vectors $\beta_1,...,\beta_n$ generate $\R^d$, because $||\cdot ||$ is a norm. Indeed, if \\ $\Span_{\R}(\beta_1,...,\beta_n)\ne \R^d$, then there would be a nonzero vector $v\in \bigcap_{i=1}^n \beta_i^{\bot}$ and hence $||v||=0$.
Let $P$ be the $n\times d$ matrix with rows equal to $\beta_1,...,\beta_n$ respectively, and let $v_1,...,v_d\in \Q^n$ be the columns of $P$. Since the rows of $P$ generate $\R^d$, the vectors $v_1,...,v_d$ are linearly independent. Thanks to the following two lemmas, we will find:
\begin{enumerate}
    \item There is a symmetric $n\times n$ matrix $\overline A$ with integral coefficients such that $\ker \overline A=\Span_{\R}(v_1,...,v_d)$;
    \item There is a closed oriented good graph manifold $M$ with base surfaces of genera $g_1,...,g_n$, whose reduced plumbing matrix is $A=\frac 1r\overline A$, for some integer $r\ne 0$. Also, the base surfaces of the pieces of $M$ can be chosen to have all orbifold Euler characteristic $\chi_i<0$.
\end{enumerate}
But $\ker A=\ker \overline A$, thus Proposition \ref{prop: nonvanishing norm} proves that $||\cdot ||$ coincides with $x^{nv}$ in the basis $(v'_1,...,v'_d)$ of $\ker A$, where the $j$-th component of $v'_i$ is $-\frac 1{\chi_j}$ times the $j$-th component of $v_i$. To construct a rational basis of $H_2(M)$ as desired, just take representatives of $v'_1,...,v'_d$ and extend them to a basis by adding enough linearly independent classes of vertical surfaces.

Finally, for the fibering statement, recall that a class $\alpha\in H_2(M)$ is fibered if and only if $\psi_i(\alpha)\ne 0$ for each $i=1,...,n$ (see the proof of Proposition \ref{prop: fibering}). In particular, $M$ does not fiber over the circle if and only if there is an $i$ so that $\psi_i(\alpha)=0$ for every $\alpha\in H_2(M)$, i.e. $\beta_i=0$. Hence, if we want to construct an $M$ fibering over the circle, we just need to repeat the construction after discarding the $\beta_i$s which are zero. Oppositely, if we want $M$ not to fiber over the circle, then it is enough to repeat the construction with $\beta_1,...,\beta_n, \beta_{n+1}=0$.
\ep

\bl For every $v_1,...,v_d\in\Q^n$ linearly independent vectors, there is an $n\times n$ integral symmetric matrix $\overline A$ such that $\ker \overline A=\Span_{\Q}(v_1,...,v_d)$.
\el
\bp Let $B=(b_{i,j})_{i,j=1,..,n}$ be the matrix defined as $$b_{i,j}=\begin{cases}
    1 \text{ if $i=j>d$}\\
    0 \text{ otherwise}.
\end{cases}$$ The kernel of $B$ is the subspace generated by the first $d$ vectors of the canonical basis of $\Q^n$. Complete $v_1,...,v_d$ to a basis $v_1,...,v_d,w_{d+1},...,w_n$ of $\Q^n$, and let $Q$ be matrix whose columns are the vectors $v_1,...,v_d,w_{d+1},...,w_n$. The matrix $\overline B:=\left(Q^{-1}\right)^tBQ^{-1}$ is a symmetric rational matrix whose kernel is exactly $\Span_{\Q}(v_1,...,v_d).$ After multiplying $\overline B$ by a sufficiently large integer, we find an integral matrix $\overline A$ with the same kernel as $\overline B$.
\ep

\bl Let $\overline A$ be a symmetric $n\times n$ integral matrix and let $g_1,...,g_n$ be nonnegative integers. There exists a closed oriented good graph manifold $M$ with base surfaces of genera $g_1,...,g_n$, such that the reduced plumbing matrix $A$ of $M$ satisfies $$\overline A=rA$$ for some integer $r\ne 0$. Also, the base surfaces of the pieces of $M$ can be chosen to have all orbifold Euler characteristic $\chi_i<0$. 
\el
\bp If $\overline A$ is the zero matrix, we can just choose $M$ as a disjoint union of manifolds of the form $\Sigma\times S^1$, where $\Sigma$ is any closed oriented surface. Hence suppose that $\overline A$ is nonzero. Let $N$ be the least common multiple of the nonzero entries of $\overline A$, and let $A$ be the matrix $\frac 1N \overline A$. 

First, we find a closed oriented good graph manifold $M'$ whose reduced plumbing matrix coincides with $A$ on the nondiagonal entries. For $i=1,...,n$, let $M_i'=\Sigma_i\times S^1$, with $\Sigma_i$ the compact oriented surface of genus $g_i$ and as many boundary components as the number of the nonnull and nondiagonal entries in the $i$-th row of $A$. Let then $M'$ be the graph manifold with pieces $M_1',...,M_n'$, where a boundary torus of $M_i'$ is glued to a boundary torus of $M_j'$ if and only if $i\ne j$ and $A_{i,j}\ne 0$. In this case, we glue the torus through any diffeomorphism such that the regular fibers on the two sides intersect algebraically $1/A_{i,j}$ times.

The diagonal entries of the reduced plumbing matrix of $M'$ are the rationals $e_1',...,e_n'$.
For $i=1,...,n$, let $\alpha_i,\beta_i$ be coprime integers such that $-\beta_i/\alpha_i=A_{i,i}-e_i'$ (to check the signs, see Definition \ref{def: invariants}). Let $M$ be the graph manifold obtained from $M$ by performing an $(\alpha_i,\beta_i)$-Dehn surgery on each $M_i'$. The manifold $M$ is the desired graph manifold.

If some $\chi_i\ge 0$, we can just execute more surgeries on $M_i'$, e.g. we can execute $m$ $(m\alpha_i,\beta_i)$-Dehn surgeries, instead of just one $(\alpha_i,\beta_i)$-Dehn surgery, for any $m\in \N$. 
\ep

\bc\label{cor: planar norms} If $||\cdot ||$ is a norm on $\R^2$ whose unit ball is a polygon with rational vertices, there exist a closed oriented good graph manifold $M$ and a rational basis on $H_2(M)$, such that the nonvanishing Thurston norm on $H_2^{nv}(M)$ coincides with $||\cdot ||$ in the induced basis. Moreover, we can ask for either all or none of the top-dimensional Thurston cones of $H_2(M)$ to be fibered.
\ec
\bp The same reasoning as the one in \cite[Chapter 4]{norm} shows that every such norm is a sum of absolute values of linear functionals $\Q^2\to\Q$. The statement follows from Theorem \ref{thm: every sum is graph} by choosing all $g_i=0$.
\ep

Unfortunately, we will see in the next section that the analogous statement does not hold in higher dimension.

\section{Linear cellular decompositions of the sphere induced by Thurston norms}\label{sec: cellular decompositions}

Now, we aim to explore geometrically the nonvanishing Thurston norms of graph manifolds. We will describe how to recover the shape of their unit polyhedra and observe a key property of them, that we will call ``completeness". Looking for completeness of the unit balls will allow us to show that the family of norms realisable as nonvanishing Thurston norms of graph manifolds is strictly contained in the family of nonvanishing Thurston norms of compact oriented $3$-manifolds (see Example \ref{exa: 3-chain}).

\bigskip

\noindent A norm on $\R^n$ with polyhedral unit ball induces a partition of $\R^n$ into convex cones where the norm equals a linear functional. As a first tool to study these norms, we forget the linear functionals and just look at the partition that the cones induce on the sphere $S^{n-1}$.

\bd[Linear cellular decompositions of $S^{n-1}$]\label{def: completeness} Let $P\subset \R^n$ be a convex (solid) polyhedron containing the origin in its interior. Let $\pi:\R^n-\{0\}\to S^{n-1}$ be the radial projection. The cellular decomposition of $\partial P$ into faces induces a cellular decomposition of $S^{n-1}$ via the projection $\pi$. If $F$ is a $k$-dimensional (open) face of $\partial P$, then $\pi(F)$ is a $k$-dimensional cell in $S^{n-1}$, contained in the intersection of $S^{n-1}$ with a $(k+1)$-dimensional subspace of $\R^n$, for every $k\le n-1$. We call such a cellular decomposition of $S^{n-1}$ a \emph{linear cellular decomposition} and we indicate it with the symbol $\mathcal F(P)$. For $k=0,...,n-1$, $\mathcal F^k(P)$ will indicate the $k$-skeleton of $\mathcal F(P)$.

Since the cellular decomposition is linear, the ($n-2$)-dimensional skeleton is contained in the intersection of $S^{n-1}$ with some hyperplanes of $\R^n$, namely there exist $\beta_1,...,\beta_k\in \R^n-\{0\}$ such that $$\mathcal F^{n-2}(P)\subset S^{n-1}\cap\left(\bigcup_{i=1}^k \beta_i^{\bot}\right).$$
If the inclusion above is an equality, we say that the linear cellular decomposition is \emph{complete}, and the polyhedron $P$ is \emph{complete}. When the $\beta_i$ are all in $\Q^n$, we refer to the linear cellular decomposition as a \emph{rational cellular decomposition}.

If $||\cdot ||$ is a norm on $\R^n$ whose unit ball is a polyhedron $P$, we call $\mathcal F(||\cdot||):=\mathcal F( P)$. Such a linear cellular decomposition is \emph{symmetric}, meaning that $\mathcal F^k(||\cdot||)$ is preserved by the antipodal map $v\to -v$, for every $k=0,...,n-1$.
\ed

\begin{figure}[h]
    \centering
    \includegraphics[width=0.49\linewidth]{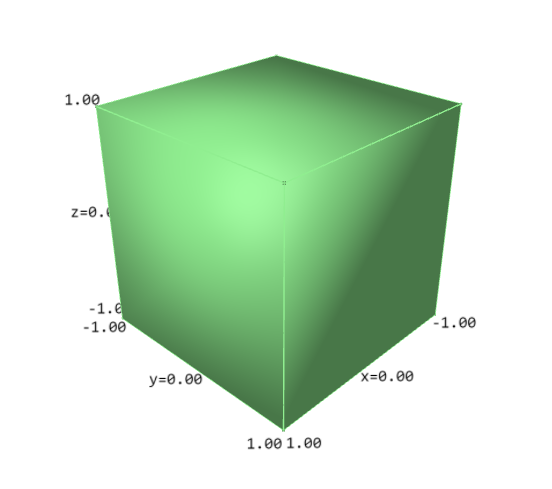}
    \includegraphics[width=0.49\textwidth]{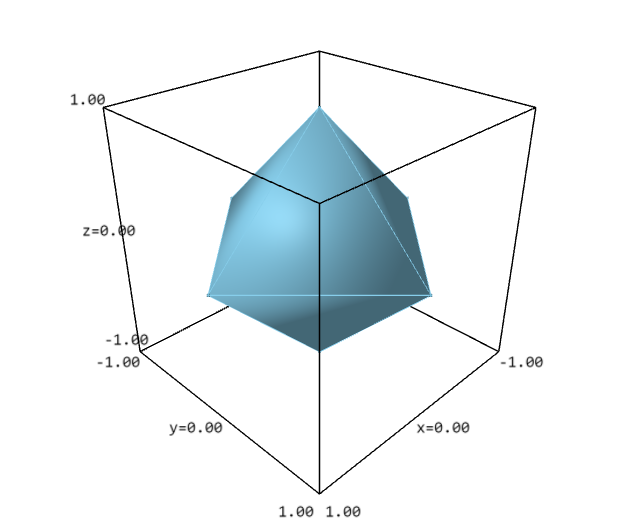}
    \caption{}
    \label{fig: cube-octa}
\end{figure}

\bexa\label{exa: cube} Let $C$ be the cube $[-1,1]^3\subset \R^3$, see Figure \ref{fig: cube-octa} left. The linear cellular decomposition of $S^2$ induced by $C$ is rational and symmetric but not complete, indeed $$\mathcal F^1(C)\subset S^2\cap(\{x\pm y=0\}\cup\{x\pm z=0\}\cup\{y\pm z=0\})$$ but the equality does not hold.

Let $D\subset \R^3$ be the octahedron with vertices $(\pm 1,0,0)$, $(0,\pm 1,0)$ and $(0,0,\pm 1)$, see Figure \ref{fig: cube-octa} right. The linear cellular decomposition of $S^2$ induced by $C$ is rational and complete, indeed $$\mathcal F^1(D)=S^2\cap( \{x=0\}\cup\{y=0\}\cup\{z=0\}).$$ 
Observe that $C$ is the unit ball of $$||(x,y,z)||=\max(|x|,|y|,|z|),$$ whilst $D$ is the unit ball of $$||(x,y,z)||=|x|+|y|+|z|.$$
\eexa

\bprop\label{prop: completeness} Let $\beta_1,...,\beta_k\in \R^n-\{0\}$ be a set of generators of $\R^n$. Let $||\cdot||$ be the norm $$||v||:=\sum_{i=1}^k|\langle \beta_i,v\rangle|.$$ Then, $||\cdot||$ has polyhedral unit ball, the linear cellular decomposition of $S^{n-1}$ induced by $||\cdot||$ is complete and \begin{equation}\label{equation: completeness}
    \mathcal F^{n-2}(||\cdot||)= S^{n-1}\cap\left(\bigcup_{i=1}^k \beta_i^{\bot}\right).
\end{equation}
\eprop
\bp If $C$ is a connected component of $\R^n-\bigcup_{i=1}^k \beta_i^{\bot}$, each $\langle \beta_i,\cdot\rangle$ has a constant sign on $C$, hence $||\cdot ||$ coincides with a linear functional on $C$. This implies that the unit ball of $||\cdot||$ is a polyhedron, that $C$ is contained in the cone over a top-dimensional face of the unit ball, and the inclusion $$\mathcal F^{n-2}(||\cdot||)\subset S^{n-1}\cap\left(\bigcup_{i=1}^k \beta_i^{\bot}\right)$$ holds.

Let $v\in \beta_i^{\bot}-\left(\bigcup_{j\ne i}\beta_j^{\bot}\right)$, we show now that $v$ does not lie in the cone over a top-dimensional face of the unit ball of $||\cdot||$. Let $||\cdot ||'$ be the seminorm $$||\cdot||':=\sum_{j\ne i} |\langle\beta_j,\cdot\rangle|.$$ The same line of reasoning as above shows that $v$ belongs to the open cone $C'$ where $||\cdot||'$ coincides with a linear function $L:\R^n\to\R$. Since $C'$ is open, there is an $\epsilon >0$ such that both $v+\epsilon\beta_i$ and $v-\epsilon\beta_i$ lie in $C'$. We can then compute $$||v+\epsilon\beta_i||+||v-\epsilon\beta_i||=L(v+\epsilon\beta_i)+|\langle \beta_i,v+\epsilon\beta_i\rangle|+L(v-\epsilon\beta_i)+|\langle \beta_i,v-\epsilon\beta_i\rangle|=2L(v)+2\epsilon \langle \beta_i,\beta_i\rangle>2L(v)=2||v||.$$ Thus, the cone $C'$ is divided by $\beta_i^{\bot}$ into two cones, where $||\cdot||$ restricts to two different linear functionals. This shows that $S^{n-1}\cap\left(\beta_i^{\bot}-\left(\bigcup_{j\ne i}\beta_j^{\bot}\right)\right)$ is a union of $(n-2)$-dimensional cells of $\mathcal F^{n-2}(||\cdot||)$. In fact, the $(n-2)$-dimensional cells of $\mathcal F^{n-2}(||\cdot||)$ correspond exactly with the union of the $S^{n-1}\cap\left(\beta_i^{\bot}-\left(\bigcup_{j\ne i}\beta_j^{\bot}\right)\right)$, for $i=1,...,k$.

To study $\mathcal F^{n-3}(||\cdot||)$ we can now restrict to study the linear cellular decompositions that $||\cdot||$ induces on every $S^{n-2}=S^{n-1}\cap\beta_i^{\bot}$. The same reasoning above can be iterated and we can prove the statement by induction on $n$, where the base step $n=2$ is proven by the reasoning above.
\ep

\begin{figure}[t]
    \centering
    \includegraphics[width=0.40\linewidth]{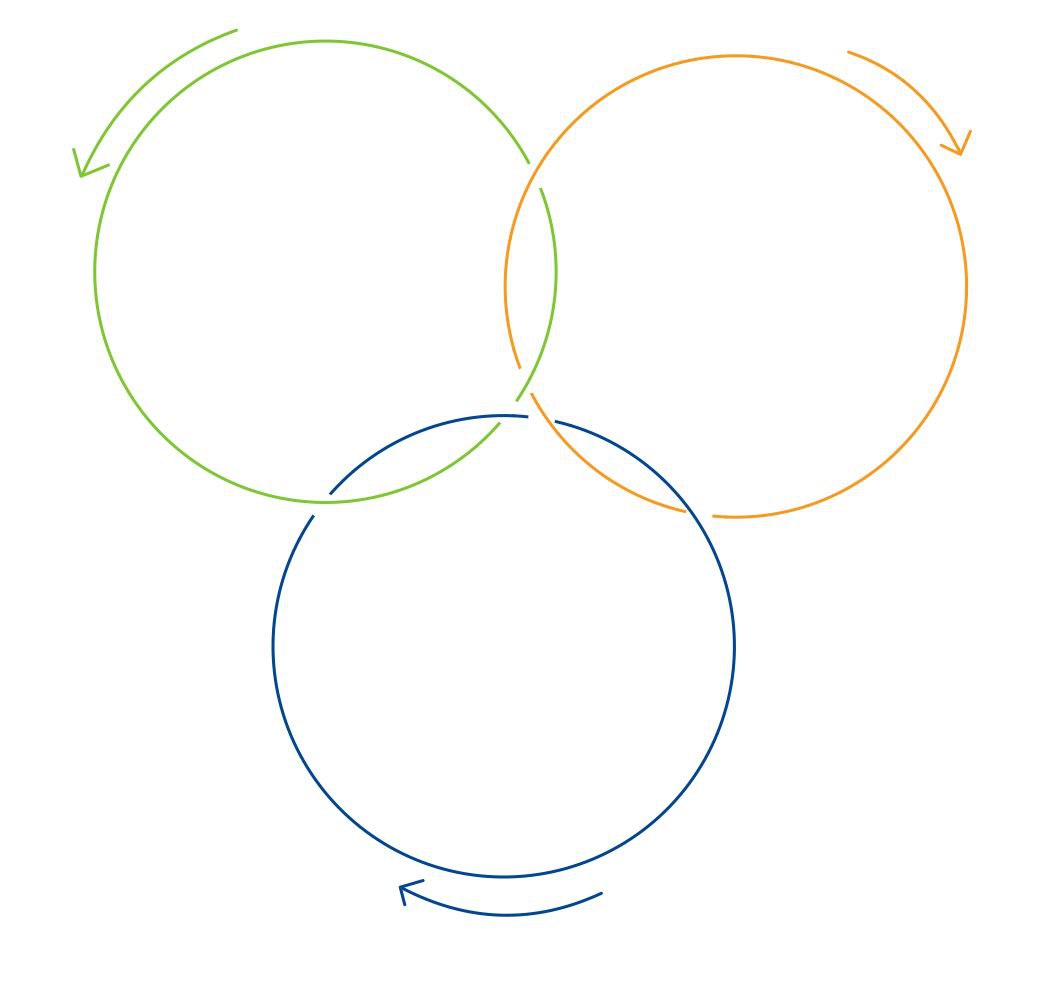}\includegraphics[width=0.50\linewidth]{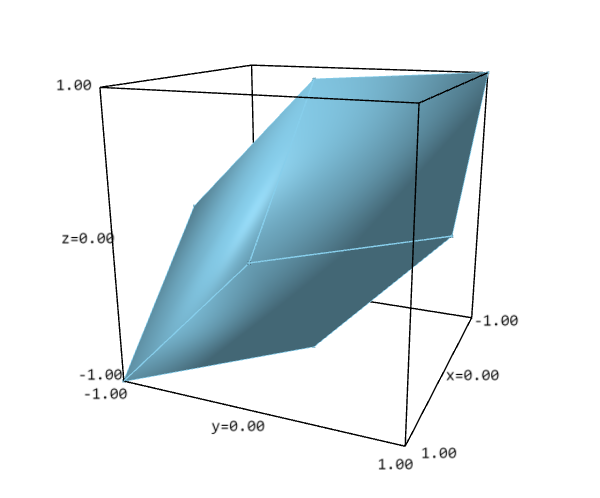}
    \caption{Left: a 3-chain link. Right: the Thurston ball of its exterior.}
    \label{fig: 3-chain}
\end{figure}

\bexa The norm $$||(x_1,...,x_n)||_{\infty}:=\max_i |x_i|$$ has unit ball $[-1,1]^n\subset\R^n$. Since $\mathcal F(||\cdot||_{\infty})$ is not complete (see Example \ref{exa: cube}), this norm cannot be represented as a sum of absolute values of linear functionals. In particular, it cannot be realised as the nonvanishing Thurston norm of a closed oriented graph manifold.
\eexa

\bexa\label{exa: 3-chain} Let $L$ be a $3$-link chain and let $M:=S^3-\overset{\circ}{N}(L)$, see Figure \ref{fig: 3-chain}. As shown in \cite{norm}, the Thurston unit ball on $H_2(M,\partial M)$ is the polyhedron with vertices $(1,0,0),(0,1,0),(0,0,1),(-1,1,1)$ and their opposites. The linear cellular decomposition of $S^2$ induced by the Thurston norm is not complete, as there are vertices with odd degree. Hence, the family of norms that are realizable as nonvanishing Thurston norm of graph manifolds is strictly included in the family of norms realizable as nonvanishing Thurston norm of $3$-manifolds.
\eexa

\br In \cite{norm}, Thurston showed that if a convex symmetric polygon in $\R^2$ has vertices $\pm v_1,...,\pm v_k$ in $\Z^2$, with $v_i\equiv v_j$ modulo $2$ for all $i,j$, then it can be realised as the dual Thurston ball of a compact orientable $3$-manifold. As pointed out by Sane in \cite{karim2}, the idea at the base of Thurston's proof is that every such polygon is the dual unit ball of an \emph{intersection norm} on the first homology group of a torus. For every intersection norm on the first homology of a torus, Thurston finds a $3$-manifold whose Thurston norm on the second homology coincides with the given intersection norm. In particular, Sane generalizes this construction to the higher-dimensional case: he considers intersection norms on the first homology of the closed orientable surface of genus $g$ and their induced polyhedra in $\R^{2g}$. As shown in \cite{karim1}, not every symmetric polyhedron is realizable this time, nor assuming the same parity condition working for polygons. 

Another way to see this is that every intersection norm can be written as a sum of absolute values of linear functionals, so its unit ball is a complete polyhedron. Hence, no noncomplete polyhedron can be realised as the unit ball of an intersection norm.
\er

Also, thanks to Theorem \ref{thm: every sum is graph} and Proposition \ref{prop: completeness}, we can state the following corollary.

\bc\label{cor: completion} Every complete symmetric rational cellular decomposition of $S^{n-1}$ is induced by the nonvanishing Thurston norm on the nonvanishing second homology of a closed oriented good graph manifold with respect to a rational basis. Also in this case, we can ask for either all or none of the top-dimensional Thurston cones to be fibered.
\ec

\begin{figure}[t]
    \centering
    \includegraphics[width=0.5\linewidth]{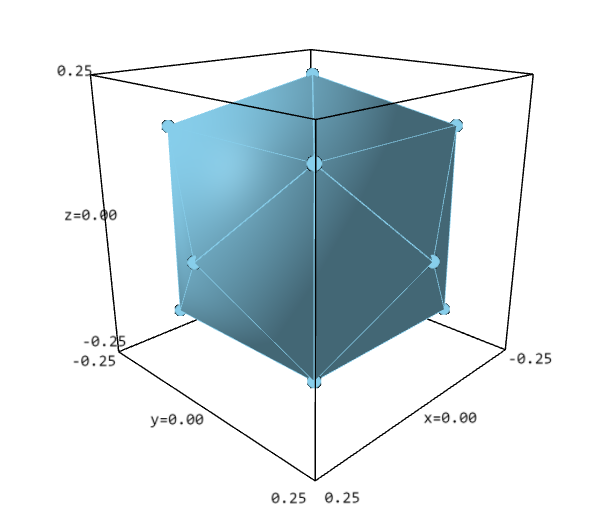}
    \caption{The polyhedron $P$, here represented for $\epsilon=\frac 14$.}
    \label{fig: nonadmissible}
\end{figure}

In particular, given a polyhedron $P$ with rational vertices in $\R^n$, we can look at the induced rational cellular decomposition of $S^{n-1}$. This cellular decomposition can be \emph{completed}: if $\beta_1,...,\beta_k\in \Q^n$ is a minimal family of vectors such that $$\mathcal F^{n-2}(P)\subset S^{n-1}\cap\left(\bigcup_{i=1}^k \beta_i^{\bot}\right)$$ then its \emph{completion} is the rational cellular decomposition whose $(n-2)$-skeleton is exactly $S^{n-1}\cap\left(\bigcup_{i=1}^k \beta_i^{\bot}\right)$. For instance, the completion of the cellular decomposition induced by the cube in Figure \ref{fig: cube-octa}-left is the cellular decomposition induced by the polyhedron in Figure \ref{fig: nonadmissible}.
The completion is induced by the nonvanishing Thurston norm of some closed oriented good graph manifold. We can summarize this discussion with the following corollary.

\bc\label{cor: completion2} Let $P\subset \R^n$ be a symmetric polyhedron with rational vertices. There is a closed oriented good graph manifold $M$, such that there is a rational basis in which $H_2^{nv}(M)\cong \R^n$ and the partition of $H_2^{nv}(M)$ into cones over the faces of the nonvanishing Thurston unit ball is a refinement of the partition of $\R^n$ into the cones over the faces of $P$. Moreover, we can ask for either all or none of the top-dimensional Thurston cones of $H_2(M)$ to be fibered.
\ec

It is worth stressing that both Corollaries \ref{cor: completion} and \ref{cor: completion2} are here established with respect to the choice of a \emph{rational basis} for the second homology of the relative graph manifold. In order to show the same results with respect to integral bases, it would be enough to show that for every $\beta_1,...,\beta_n$ integral generators of $\Z^d$ there exist a closed oriented good graph manifold $M$, an integral basis of $H_2(M)$ and a natural number $N$, such that the nonvanishing Thurston norm is $N\sum_1^n |\langle\beta_i,\cdot\rangle|$ in the induced integral basis. Pursuing this strategy would demand a closer look to the relationship between $\psi(H_2^{nv}(M;\Z))$ and $\ker A\cap \Z^n$ (see Proposition \ref{prop: computation}).

\begin{figure}[t]
        \centering
        \includegraphics[width=0.5\linewidth]{octahedron.png}\includegraphics[width=0.42\linewidth]{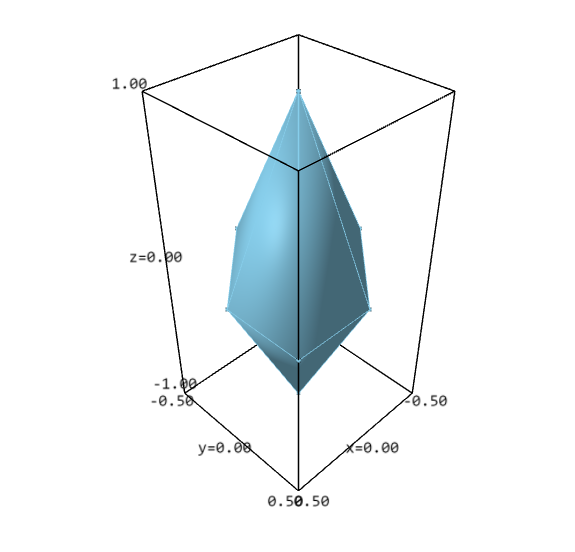}
        \caption{Left: the octahedron $P$, unit ball of $|x|+|y|+|z|$. Right: the polyhedron $P'$, unit ball of $|x|+|y|+|z|+|x+y|$. The intersections of $P$ and $P'$ with $\{x+y=0\}$ coincide. Then, the rest of $P'$ is obtained by ``deflating" $P$ symmetrically with respect to $\{x+y=0\}$.}
        \label{fig: deflating}
    \end{figure}

\begin{construction}\label{construction}
    We will now give an algorithm to visualize the unit ball of a sum of absolute values of linear functionals. Let $$||\cdot||:=\sum_{i=1}^k|\langle \beta_i,\cdot\rangle|,$$ where $\beta_1,...,\beta_k\in\R^n$ generate $\R^n$ (this corresponds to $||\cdot||$ being a norm). The algorithm is by induction on $k$.

    Suppose $k=n$, meaning that the vectors $\beta_1,...,\beta_n$ are linearly independent. For each $i$, \\ $\Span_{\R}\{\beta_j\,|\; j\ne i\}$ has dimension $n-1$, then choose $\overline \beta_i$ to be a non null vector in $\bigcap_{j\ne i}\beta_j^{\bot}$, normalized to have $||\overline\beta_i||^2=|\langle \overline\beta_i,\overline\beta_i\rangle|=1$. The unit ball of $||\cdot||$ is the polyhedron with vertices $\pm\overline \beta_1,...,\pm\overline\beta_n$ in $\R^n$. 

    Suppose that $k>n$. After possibly relabeling, suppose $\beta_k\in\Span_{\R}(\beta_1,...,\beta_{k-1})$. Write $$||\cdot||=||\cdot||'+|\langle \beta_k,\cdot\rangle|$$ and suppose by induction that we can describe the unit ball $P'$ of $||\cdot||'$, $\pm v_1',...,\pm v_l'$ being its vertices. The unit ball $P$ of $||\cdot ||$ is obtained from $P'$ by a slight ``symmetrical deflation" of $P'$ with respect to $\beta_k^{\bot}$, but letting $P'\cap \beta_k^{\bot}$ remain unaffected, see Figure \ref{fig: deflating}. Indeed, the $1$-codimensional polyhedron $P'\cap \beta_k^{\bot}$ coincides with $P\cap \beta_k^{\bot}$, whilst each vertex $\pm v_i'$ of $P'$ gives rise to a vertex $\pm v_i$ of $P$ satisfying $$v_i:=\frac 1{1+|\langle \beta_k,v_i'\rangle|}v_i'.$$
    Apart from $\pm v_1,...,\pm v_l$, the other vertices of $P$ are the vertices of $P'\cap \beta_k^{\bot}$.
\end{construction}

\br Construction \ref{construction} makes us compute the unit ball of the nonvanishing Thurston norm of a closed oriented good graph manifold $M$, thanks to Proposition \ref{prop: nonvanishing norm}. Then, by Proposition \ref{prop: dimension}, we can recover the full Thurston polyhedron on $H_2(M)$. 

When $M$ is a closed oriented graph manifold, but not good, it admits a good double cover $\widetilde M$. The pull-back map $\pi^*: H^1(M)\to H^1(\widetilde M)$ is injective and for every $\alpha\in H^1(M)\cong H_2(M)$ we have \begin{equation}\label{equation: cover}
    x_{\widetilde M}(\pi^*\alpha)=2x_M(\alpha).
\end{equation} Identity (\ref{equation: cover}) was proved by Gabai in \cite{gabai1}. In particular, the Thurston ball of $H_2(M)$ can be fully recovered by intersecting the Thurston ball of $H_2(\widetilde M)$ with the subspace $\pi^*H_2(M)$ and rescaling it by a factor $2$.
\er

\br\label{rem: icosaedro} By Lemma \ref{lemma: norm} and Proposition \ref{prop: completeness}, the nonvanishing Thurston unit ball on the nonvanishing second homology of a closed oriented graph manifold is a symmetric complete polyhedron, and it has rational vertices. Unfortunately, not every symmetric complete polyhedron with rational vertices can be obtained as the nonvanishing Thurston ball of some closed oriented graph manifold. As a counterexample, consider the complete polyhedron $P\subset\R^3$ in Figure \ref{fig: nonadmissible} with vertices $$\pm\left(\frac 14,0,0\right),\pm\left(0,\frac 14,0\right),\pm\left(0,0,\epsilon\right),\pm\left(\frac 16,\frac 16,\frac 16\right),\pm\left(\frac16,\frac16,-\frac16\right),\pm\left(\frac 16,-\frac16,\frac16\right),\pm\left(\frac16,-\frac 16,-\frac16\right)$$ for some $\epsilon>0$ to be determined. Suppose that $P$ is the unit ball of a norm $||\cdot||$ that can be expressed as a sum of absolute values of linear functionals. Then, thanks to equality (\ref{equation: completeness}), there are rational numbers $a,b,c,d,e,f$ such that for every $x,y,z\in\R$ $$||(x,y,z)||=a|x+y|+b|x+z|+c|y+z|+d|x-y|+e|x-z|+f|y-z|.$$ By imposing each vertex apart from $\pm(0,0,\epsilon)$ to have unit norm, we get that necessarily $a=b=c=d=e=f=1$. But the unit ball of $||\cdot||$ for such choice of coefficients gives the vertices $\pm\left(0,0,\frac 14\right)$. In particular, there are infinite values of $\epsilon$ such that the polyhedron $P$ is not the unit ball of a sum of absolute values of linear functionals.
\er

\end{document}